\documentclass[oneside]{amsart}
\newcommand{\formatswitch}{preprint}

\usepackage{amssymb}
\usepackage{ifthen}
\usepackage{verbatim}
\usepackage{psfig}
\usepackage[all]{xy}

\newcommand{\tref}[1]{(\ref{#1})}

\ifthenelse{\equal{\formatswitch}{big}}{
\setlength{\textwidth}{432pt}
\setlength{\oddsidemargin}{18.8775pt}

\setcounter{tocdepth}{1}
}{
\ifthenelse{\equal{\formatswitch}{paper}}{

\setcounter{tocdepth}{1}
}{
\setcounter{tocdepth}{2}
}
}
\DeclareMathAlphabet\EuScript{U}{eus}{m}{n}
\DeclareMathAlphabet\EuScriptb{U}{eus}{b}{n}

\newcommand{\sscr}[1]{\EuScript{#1}}

\newcommand{\claimenum}{\renewcommand{\theenumi}{\alph{enumi}}
 \renewcommand{\labelenumi}{\textit{(\theenumi)}}
 \renewcommand{\theenumii}{\roman{enumii}}
 \renewcommand{\labelenumii}{\textit{(\theenumii)}}
 \begin{enumerate}}
\newcommand{\claimenumend}{\end{enumerate}}

\newcommand{\romanenum}{\renewcommand{\theenumi}{\roman{enumi}}
 \renewcommand{\labelenumi}{\textit{(\theenumi)}}
 \renewcommand{\theenumii}{\alph{enumii}}
 \renewcommand{\labelenumii}{\textit{(\theenumii)}}
 \begin{enumerate}}
\newcommand{\romanenumend}{\end{enumerate}}

\newtheorem{dummy}{realdumb}
\newtheorem{thm}{Theorem}

\newtheorem{lemma}[dummy]{Lemma}

{\theoremstyle{definition} }

{\theoremstyle{definition} }

\renewcommand{\text}{\mathrm}

\newcommand{\strutdepth}{\dp\strutbox}
\newcommand{\marginalnote}[1]
   {\strut\vadjust{\kern-\strutdepth\domarginalnote{#1}}}
\newcommand{\domarginalnote}[1]{\vtop to \strutdepth{
  \baselineskip\strutdepth
   \vss\llap{ #1\ \ }\null}}  
\newcounter{showlabelflag}
\newcounter{makelabelflag}
\newcommand{\showlabels}{\setcounter{showlabelflag}{1}}

\newcommand{\makelabels}{\setcounter{makelabelflag}{1}}
\newcommand{\hidelabels}{\setcounter{showlabelflag}{2}}
\newcommand{\mylabel}[1]{
  \ifthenelse{\value{makelabelflag}=1}
    {\label{#1}}{}
  \ifthenelse{\value{showlabelflag}=1}
    {\marginpar{#1}}{}\relax}

\ifthenelse{\equal{\formatswitch}{draft}}{
\showlabels
}{\relax}

\newcommand{\scr}{\sscr}

\newcommand{\mymargin}[1]{
  \ifthenelse{\value{showlabelflag}=1}
    {\marginpar{#1}}{}\relax}

\newcounter{enumo}

\newcommand{\RRsh}{\kern -1 pt \Rsh}

\newcounter{keepitemnum}

\newcounter{keepitemnumm}

\begin{document}

\bibliographystyle{amsplain}
\title[Coherence of Associativity]{Coherence of Associativity in 
\\ Categories with Multiplication
}
\thanks{AMS
Classification (2000): primary 18D10, secondary 20F05.}
\author{MATTHEW G. BRIN}
\date{September 4, 2004}

\CompileMatrices


\makelabels
\hidelabels
\maketitle


\section{Introduction}\mylabel{IntroSec}

To say that \(\scr{C}\) is a category with (functoral)
multiplication means that there is a functor
\(\otimes:\scr{C}^2\rightarrow \scr{C}\) called the multiplication
where \(\scr{C}^2\) is the category of pairs of objects and pairs of
morphisms from \(\scr{C}\).  [More technically, \(\scr{C}^2\) is the
category of functors and natural transformations from \(2\) to
\(\scr{C}\) where \(2\) is the category with objects 0 and 1, and
the only morphisms are the identity morphisms.]  Examples of
functoral multiplications are cross products, tensor products, free
products and so forth on those categories where those products
exist.

For most examples it is rarely the case that
\mymargin{AssocLaw}\begin{equation}\label{AssocLaw}
A\otimes(B\otimes C) = (A\otimes B)\otimes C \end{equation} is
literally true, and what is usually the case is that there is a
natural isomorphism \(\alpha\) from the functor
\(F:\scr{C}^3\rightarrow \scr{C}\) defined by \(F(A,B,C)=A\otimes
(B\otimes C)\) to the functor \(G:\scr{C}^3\rightarrow \scr{C}\)
defined by \(G(A,B,C)=(A\otimes B)\otimes C\).  In the most common
cases, there is an obvious candidate for natural isomorphism
\(\alpha\) and it is a triviality to define.

The usual statement that ``all associativity laws follow from the
associativity law given in \tref{AssocLaw}'' translates into a claim
that if \(H\) and \(K\) are two functors from \(\scr{C}^n\) to
\(\scr{C}\) that are built by combining \(n\) variables in the same
order with \(n-1\) applications of \(\otimes\) and that differ only
in the pattern of parentheses, then there is a natural isomorphism
from \(H\) to \(K\) that is derivable in some sensible way from
\(\alpha\).  The problem might be that there is more than one way to
build such an isomorphism from \(\alpha\), raising the possibility
that different ways will result in different isomorphisms.

This problem was first considered by MacLane in
\cite{MacLane:coherence} where he defined the condition {\itshape
coherence} of such an \(\alpha\) to mean that any two expressions
built from \(\otimes\) using the same variables in the same order
and differing only in the distribution of parentheses are connected
by a unique natural isomorphism derivable from \(\alpha\) using a
prescribed set of constructions.  In \cite{MacLane:coherence} it is
proven that coherence is achieved from the naturality of \(\alpha\)
and one hypothese that a certain (now famous) pentagonal diagram
commutes.

The purpose of this paper is to show that the hypothesis that the
pentagonal diagram commute can be dispensed with if the prescribed
set of constructions for building natural isomorphisms from
\(\alpha\) is restricted.  Thus we do not prove a strengthening or
generalization of MacLane's theorem.  It is simply a different
theorem.

Beyond the statement and proof of this theorem, the paper has a
second purpose which is to point out the connection between
MacLane's theorem on coherence and combinatorial group theory.  This
is discussed in the last section where we point out that MacLane's
theorem can be viewed as giving a presentation of a certain group in
terms of generators and relations.

\section{Statement}\label{StatementSec}

The constructions in \cite{MacLane:coherence} for building
isormophisms from \(\alpha\) are extremely natural.  (Overuse of the
word {\itshape natural} here is unavoidable.)  The restrictions on
the constructions in this paper lack a certain symmetry.  Thus our
result suffers from a certain aesthetic inferiority.  We now give
some details and start with some preliminary technicalities.

If \(\beta\) is a natural transformation from a functor
\(F:\scr{A}\rightarrow \scr{B}\) to a functor \(G:\scr{A}\rightarrow
\scr{B}\), then we can view \(\beta\) as a functor from \(\scr{A}\)
to \(\scr{B}^{\overline{2}}\), the category of functors from
\(\overline{2}\) to \(\scr{B}\) in which \(\overline{2}\) is the
category with objects 0 and 1 and only one non-identity morphism
that goes from 0 to 1.  The category \(\overline{2}\) is just the
category whose objects are 0 and 1 and whose morphisms correspond to
the partial order \(\le\); while \(\scr{B}^{\overline{2}}\) is just
the category whose objects are the morphisms of \(\scr{B}\) and
whose morphisms are the commutative squares in \(\scr{B}\).  If
\(S\) is the ``source'' functor from \(\scr{B}^{\overline{2}}\) to
\(\scr{B}\) in which \(S(f:X\rightarrow Y)=X\) and \(T\) is the
``target'' functor in which \(T(f:X\rightarrow Y)=Y\), then
\(S\beta=F\) and \(T\beta=G\).

Any functor \(F:\scr{A}\rightarrow \scr{B}\) induces a functor
\(F^{\overline{2}}:\scr{A}^{\overline{2}}\rightarrow
\scr{B}^{\overline{2}}\).

In \cite{MacLane:coherence} isomorphisms are built from 
\(\alpha: A\otimes(B\otimes C) \rightarrow (A\otimes B)\otimes C \)
by four processes.  The one that we will restrict is as follows.

If \(\beta\) is a natural transformation from functor \(F\) to
functor \(G\) that each go from \(\scr{C}^m\) to \(\scr{C}\) and
\(\gamma\) is a natural transformation from \(H\) to \(K\) that each
go from \(\scr{C}^n\) to \(\scr{C}\), then we can form
\(\beta\otimes \gamma\) going from \(F\otimes H\) to \(G\otimes K\)
by composing \[\beta\times \gamma:\scr{C}^m\times \scr{C}^n
\rightarrow \scr{C}^{\overline{2}} \times \scr{C}^{\overline{2}}\]
with \[\otimes^{\overline{2}}: \scr{C}^{\overline{2}} \times
\scr{C}^{\overline{2}} \rightarrow \scr{C}^{\overline{2}}.\]

The operation \(\otimes\) on transformations can be used for the
following.  Let \(\mathbf{1}\) denote the identity transformation
from the identity functor on \(\scr{C}\) to itself.  We then can
form \(\alpha\otimes \mathbf{1}\),
\((\alpha\otimes\mathbf{1})\otimes \mathbf{1}\) and so forth where,
inductively, \(\alpha_0=\alpha\) and \(\alpha_i=\alpha_{i-1}\otimes
\mathbf{1}\).  Thus \[\alpha_1 = \alpha\otimes \mathbf{1} :(A\otimes
(B\otimes C))\otimes D \rightarrow ((A\otimes B)\otimes C)\otimes
D\] with similar descriptions of other \(\alpha_i\).  We can refer
to \(\beta\otimes \mathbf{1}\) as the {\itshape right stabilization}
of \(\beta\).  We refer to the \(\alpha_i\) as the iterated right
stabilizations of \(\alpha\).

The assumptions in \cite{MacLane:coherence} are that the
transformations form a category closed (among other things) under
the operation \(\otimes\) on transformations.  In this paper, we
will only make use of the operation \(\otimes\) on transformations
to create right stabilizations.  All other constructions from
\cite{MacLane:coherence} will be used here.  We now go on to the
others.

From \mymargin{AlphaBase}\begin{equation}\label{AlphaBase}
\alpha:A\otimes(B\otimes C)\rightarrow (A\otimes B)\otimes C
\end{equation} we can create \[\alpha':(A\otimes B)\otimes(C\otimes
D) \rightarrow ((A\otimes B)\otimes C)\otimes D\] from
\tref{AlphaBase} by replacing \(A\) in \tref{AlphaBase} by the
product of two variables.  Similarly, we get \[\alpha'':A\otimes
((B\otimes C)\otimes D) \rightarrow (A\otimes (B\otimes C))\otimes
D\] by replacing \(B\) in \tref{AlphaBase} by the product of two
variables.  These are both examples of {\itshape instances} of
\(\alpha\).  More generally, we can replace any variable in
\tref{AlphaBase} on both sides by identical expressions involving
\(\otimes\).

Technically, an {\itshape instance} of a transformation is created
by precomposing the transformation with a functor.  Now if
\(\beta:\scr{A}\rightarrow \scr{B}^{\overline{2}}\) is a natural
transformation from \(F=S\beta\) to \(G=T\beta\), and if
\(H:\scr{D}\rightarrow\scr{A}\) is any functor, then \(\beta D\) is
a natural transformation from \(FD=S\beta D\) to \(GD=T\beta D\) and
can be viewed as an instance of \(\beta\).

In our setting, we will take instances of the iterated right
stabilizations \(\alpha_i\) of \(\alpha\).  The iterated
stabilization of \(\alpha_i\) connects functors from \(\scr{C}^{n}\) to
\(\scr{C}\) where \(n=i+3\).  Instances can be created by
precomposing the stabilizations with compositions of functors such
as \[(X_1,X_2, \ldots, X_j, X_{j+1}, \ldots, X_{m+1}) \mapsto (X_1,
X_2, \ldots, X_j\otimes X_{j+1}, \ldots, X_{m+1})\] from
\(\scr{C}^{m+1}\) to \(\scr{C}^m\) for various values of \(m\) and
\(j\).

We will also postcompose a transformation with a functor.  If 
\(\beta:\scr{A}\rightarrow \scr{B}^{\overline{2}}\) is a natural
transformation and \(J:\scr{B}\rightarrow\scr{E}\) is a functor,
then \(J\beta\) represents the composition of \(\beta\) with
\(J^{\overline{2}}\).  This construction can yield an instance (of
another transformation) by accident, and we will exploit this.

Another operation for constructing isomorphisms from \(\alpha\) is
that of composition.  If \(F\), \(G\) and \(H\) are all functors
from \(\scr{A}\) to \(\scr{B}\), if \(\beta\) is a natural
transformation from \(F\) to \(G\) and \(\gamma\) is a natural
transformation from \(G\) to \(H\), then there is an obvious
composition \(\gamma\beta\) that is a natural transformation from
\(F\) to \(H\).  Composition commutes with right stabilization.

The final operation for constructing isomorphisms from \(\alpha\) is
that of inversion.  Since \(\alpha\), its stabilizations and its
instances are all isomorphisms, they are all invertible.  Note that
inversion commutes with instance and stabilization and behaves in
the usual way with respect to composition:
\((\beta\gamma)^{-1}=\gamma^{-1}\beta^{-1}\).  

We can now state our result.

\begin{thm}\label{MainThm} Let \(\scr{C}\) be a category with
functoral multiplication \(\otimes : \scr{C}^2\rightarrow \scr{C}\).
Let \(\alpha\) be a natural isomorphism from \(A\otimes (B\otimes
C)\) to \((A\otimes B)\otimes C\).  If \(E\) and \(F\) are two
expressions in \(n-1\) appearances of \(\otimes\) and \(n\)
different variables in the same order that differ only in the
arrangement of parentheses, then there exists a unique natural
isomorphism constructable from \(\alpha\) as a composition of
instances of iterations of right stabilizations of \(\alpha\) and
\(\alpha^{-1}\).  \end{thm}

\section{Proof}\mylabel{ProofSec}

The proof of Theorem \ref{MainThm} is essentially the proof of
Theorem 3.1 of \cite{MacLane:coherence} with more attention paid to
some details.  We will include the entire proof since a set of
instructions on modifying the proof in \cite{MacLane:coherence}
would be unreadable.

We will discuss expressions endlessly.  For us an {\itshape
expression} in \(n\) variables is a fully parenthesized alternation
of the variables \(X_1\), \dots, \(X_n\) in that order with \(n-1\)
appearances of the operation \(\otimes\).  Inductively, the variable
\(X_1\) is the only expression in 1 variable, and if \(F\) and \(G\)
are expressions in \(m\) and \(n\) variables, respectively, then
\((F\otimes \overline{G})\) is an expression in \(m+n\) variables
where \(\overline{G}\) is the expression \(G\) with all the
subscripts of its variables raised uniformly by \(m\).  We will omit
the bar from the second expression from now on since the meaning
will always be clear.

An expression is {\itshape trivial} if it has only one variable.  We
reserve the symbol \(\mathbf{I}\) to symbolize the trivial
expression \(X_1\).

A non-trivial expression \(E\) breaks uniquely as \((F\otimes
G)\).  We say that \(E\) is {\itshape semi-normalized} if
\(E=(F\otimes \mathbf{I})\).  We can refer to \((F\otimes
\mathbf{I})\) as the right stabilization of \(F\).  Right
stabilization can be iterated and we define \((F\otimes_i
\mathbf{I})\) inductively by \((F\otimes_0
\mathbf{I})=F\) and \((F\otimes_i
\mathbf{I})=((F\otimes_{i-1}\mathbf{I})\otimes\mathbf{I})\).

An expression in \(n\) variables is {\itshape fully normalized} if
it is of the form \((\mathbf{I}\otimes_{n-1}\mathbf{I})\).  There is
only one fully normalized expression on \(n\) variables for each
\(n\) and we will denote it by \(\mathbf{I}_n\).  We have
\[\mathbf{I}_n = (\cdots (((\mathbf{I} \otimes \mathbf{I})\otimes
\mathbf{I})\otimes \mathbf{I}) \otimes \cdots \otimes \mathbf{I}
), \qquad(\hbox{\(n\) apearances of}\,\,\,\mathbf{I}).\]

If an expression \(E\) is not fully normalized, then it is uniquely
expressible as \((N\otimes_i\mathbf{I})\) where \(N\) is not
semi-normalized.  If \(E\) is not semi-normalized, then \(i=0\).
The value of \(i\) is the {\itshape normalization level} of \(E\).
Note further that \(N= (F\otimes G)\) for some \(F\) and \(G\) with
\(G\ne \mathbf{I}\).  The {\itshape weight} of \(E\) is the number
of variables used in \(G\).  

If an expression \(E\) on \(n\) variables is not fully normalized,
then its normalization level is strictly less than \(n\) and its
weight is strictly greater than 1.  We extend the definitions to say
that the normalization level of \(\mathbf{I}_n\) is \(n\) and that
its weight is 1.  (There is only one expression on two variables and
it is fully normalized, so the normalization level of an \(n\)
variable expression with \(n\ge 2\) is never \(n-1\).)

The point of all this bookkeeping is the list of observations below.
They are verified by inspecting the form of the various
\(\alpha_i\).  We say that a natural transformation from a functor
\(F\) to a functor \(G\) has \(F\) as its source and \(G\) as its
target.  We are treating expressions formally, but they represent
functors.  Thus we can talk about instances of the \(\alpha_i\) as
having expressions for source and target.  The number of variables
of the source and target of a given instance of an \(\alpha_i\) will
be the same.  In reading the following, note that we carefully
distinguish between \(\alpha_i\) and \(\alpha_i^{-1}\) and the fact
that \(\alpha_i^{-1}\) is never mentioned is significant.

{
\renewcommand{\theenumi}{A\arabic{enumi}}
\begin{enumerate}

\item\mylabel{AppUniq} If \(E\) is an \(n\) variable expression,
then for each \(i\) there is at most one instance of \(\alpha_i\)
that can have \(E\) as source.

\item\mylabel{AppULim} If \(E\) is an \(n\) variable expression,
then an instance of \(\alpha_i\) can have \(E\) as source only if
\(i\le n-3\).

\item\mylabel{AppLLim} If \(E\) is an \(n\) variable expression with
normalization level \(k\), then an instance of \(\alpha_i\) can have
\(E\) as source only if \(i\ge k\).

\item\mylabel{AppToBig} If \(E\) is an \(n\) variable expression
with normalization level \(k\) and weight \(w>1\), then an instance
of \(\alpha_i\) having \(E\) as source with \(i>k\) has a target
with normalization level \(k\) and weight \(w\).

\item\mylabel{AppJustRight} If \(E\) is an \(n\) variable expression
with normalization level \(k\) and weight \(w>1\), then there is an
instance of \(\alpha_k\) having \(E\) as source.  Further the target
of this instance of \(\alpha_k\) either has normalization level that
is greater than \(k\) or has normalization level equal to \(k\) and
weight less than \(w\).

\end{enumerate}
}

If \(E\) is an expression and a string
\(\alpha_{i_1}\alpha_{i_2}\cdots \alpha_{i_s}\) has the property
that an instance of \(\alpha_{i_s}\) has \(E\) as a source and
target \(F_s\), and for each \(j<s\) an instance of \(\alpha_{i_j}\)
has \(F_{j+1}\) as a source and target \(F_j\), then we say that the
string is a word in the \(\alpha_i\) that defines a path from \(E\)
to \(F_1\).  Note that the information in the string does not
specify which instances are used, but this is not necessary because
of \tref{AppUniq}.

If an instance of \(\alpha_i\) has source \(E\) and target \(F\),
then the instance is a natural isomorphism from the functor
represented by \(E\) to the functor represented by \(F\).  Thus in
the previous paragraph, the word in the \(\alpha_i\) defines an
isomorphism from \(E\) to \(F_1\).

It is now an easy inductive exercise to prove the following from
\tref{AppUniq}--\tref{AppJustRight}.

\begin{lemma}\mylabel{ExistAndUniqLem}Given an expression \(E\) in
\(n\) variables that is not fully normalized, then there is a unique
word \[w=\alpha_{i_1} \alpha_{i_2} \cdots \alpha_{i_s}\] satisfying
\(i_j\ge i_k\) if \(j<k\) so that \(w\) is an isomorphism from \(E\)
to \(\mathbf{I}_n\).  \end{lemma}

This proves the existence part of Theorem \ref{MainThm} since any
two expressions in \(n\) variables can be connected to
\(\mathbf{I}_n\) by an isomorphism.

We now continue with our reading of the proof from
\cite{MacLane:coherence}.  If \(F\) and \(G\) are two expressions in
\(n\) variables, then we must show that any two ``paths'' from \(F\)
to \(G\) where each step in each path is an instance of some
\(\alpha_i\) or \(\alpha_i^{-1}\), then the two paths compose to the
same natural isomorphism from \(F\) to \(G\).

Our argument begins as it does in \cite{MacLane:coherence}.  We
direct each step in the path by declaring that each step in the path
goes from the source expression of the instance of \(\alpha_i\) to
its target expression.  Take an arbitrary path from \(F\) to \(G\),
and join each each vertex \(F_i\) in the path to \(\mathbf{I}_n\) by
the path obtained from Lemma \ref{ExistAndUniqLem}.  Note that this
path is directed from \(F_i\) to \(\mathbf{I}_n\).  This creates a
diagram of which the following sample is typical.
\[
\xymatrix{
{F}\ar[r] \ar[d]_{p} & {F_1} \ar[d] & {F_2} \ar[l] \ar[d] \ar[r] 
&{F_3} \ar[d] &{G} \ar[d]^{q} \ar[l] \\
{\mathbf{I}_n} \ar@{=}[r]
&{\mathbf{I}_n}\ar@{=}[r]
&{\mathbf{I}_n}\ar@{=}[r]
&{\mathbf{I}_n}\ar@{=}[r]
&{\mathbf{I}_n}
}
\]

If it is shown that the above diagram commutes, then the path along
the top from \(F\) to \(G\) gives the same isomorphism as
\(q^{-1}p\) and the proof of Theorem \ref{MainThm} will be complete.
Thus it suffices to prove the commutativity of a single rectangle of
the form
\mymargin{TestRec}\begin{equation}\label{TestRec}
\begin{split}
\xymatrix{
{F}\ar[r]^{\alpha_j} \ar[d]_{p} &{G} \ar[d]^{q} \\
{\mathbf{I}_n}\ar@{=}[r]
&{\mathbf{I}_n}
}
\end{split}
\end{equation}
in which the top arrow is an instance of some \(\alpha_j\) and the paths
\(p\) and \(q\) are obtained from Lemma \ref{ExistAndUniqLem}.

The expression \(F\) has a normalization level \(i\) so we know that
the first step in \(p\) (the last letter expressing \(p\) as a word)
is \(\alpha_i\) and we know \(j\ge i\).  If \(j=i\), then the
uniqueness gotten from Lemma \ref{ExistAndUniqLem} says that \(p\)
and \(q\alpha_j\) are identical as words and the rectangle
\tref{TestRec} commutes.  Thus we are left with the case \(j>i\).

If \(j>i\), then the normalization level of \(G\) is also \(i\) by
\tref{AppToBig} and \(q\) also has \(\alpha_i\) as its first step.
Thus we will be done by induction on the length of \(p\) when we
show that the following rectangle commutes whenever \(j>i\).
\mymargin{SmallRec}\begin{equation}\label{SmallRec}
\begin{split}
\xymatrix{
{F}\ar[r]^{\alpha_j} \ar[d]_{\alpha_i} &{G} \ar[d]^{\alpha_i} \\
{F_1}\ar[r]^{\alpha_{j+1}}
&{G_1}
}
\end{split}
\end{equation}

The expression \(F\) equals \((N\otimes_i\mathbf{I})\) with
\(N=(H\otimes(K\otimes L))\), and \(F_1=(N'\otimes_i\mathbf{I})\)
with \(N'=((H\otimes K)\otimes L)\).  Since \(j>i\), we know that
\(G=(N''\otimes_i\mathbf{I})\) where \(N''= (H'\otimes (K\otimes
L))\) is the result of applying \(\alpha_{j-i}\) to \(N\) and so
\(H'\) is the result of applying \(\alpha_{j-i-1}\) to \(H\).
Expanding what we know about \(F_1\) and applying \(\alpha_i\) to
\(G\) gives \[\begin{split} F_1 &= (((H\otimes K)\otimes
L)\otimes_i\mathbf{I}), \qquad\mathrm{and} \\ G_1 &= (((H'\otimes
K)\otimes L)\otimes_i\mathbf{I}). \end{split}\]  This says that
\(G_1\) is the result of applying \(\alpha_{j+1}\) to \(F_1\).

This does not make \tref{SmallRec} commute.  It only says that
source and targets make sense.  That commutativity follows from the
naturality of \(\alpha_i\) can be seen by filling in the details of
\tref{SmallRec} to give the following.

\mymargin{BigRec}\begin{equation}\label{BigRec}
\begin{split}
\xymatrix{
{((H\otimes(K\otimes L))\otimes_i\mathbf{I})}
\ar[rrr]^{\alpha_j} \ar[d]_{\alpha_i} 
&&&{((H'\otimes (K\otimes
L))\otimes_i\mathbf{I})} \ar[d]^{\alpha_i} \\
{(((H\otimes K)\otimes L)\otimes_i\mathbf{I})}
\ar[rrr]^{\alpha_{j+1}}
&&&{(((H'\otimes K)\otimes
L)\otimes_i\mathbf{I})}
}
\end{split}
\end{equation}

If we define functors \(R\) and \(L\) by \[\begin{split} R(-)
&=((-\otimes(K\otimes L))\otimes_i\mathbf{I}), \qquad\mathrm{and} \\
L(-) &=(((-\otimes K)\otimes L)\otimes_i\mathbf{I}), \end{split}\]
then the specific instances of \(\alpha_j\) and \(\alpha_{j+1}\) in
\tref{BigRec} are seen to be \(R\alpha_{j-i-1}\) and
\(L\alpha_{i-j-1}\), respectively.  Further, both appearances of
\(\alpha_i\) are instances of a single natural isomorphism
\(\overline{\alpha}_i\) with source \(R\) and target \(L\), and
\(\overline{\alpha}_i\) is an instance of \(\alpha_i\).  The diagram
\tref{BigRec} commutes since the following diagram commutes by the
naturality of \(\overline{\alpha}_i\).
\[
\xymatrix{
{RH} \ar[d]_{\overline{\alpha}_i} \ar[rr]^{R\alpha_{i-j-1}} &&
{RH'} \ar[d]^{\overline{\alpha}_i} \\ 
{LH}\ar[rr]^{L\alpha_{i-j-1}} && {LH'}
}
\]
This completes the proof of Theorem \ref{MainThm}

\section{The origins of Theorem \ref{MainThm} and its proof}

Theorem \ref{MainThm} is a thinly disguised translation of the well
known fact that a certain group has a certain presentation.  The
proof that we give of Theorem \ref{MainThm} contains much of the
work from the standard (and well known) proofs of this well known
fact.

There is a group commonly known as Thompson's group \(F\) (see
\cite{CFP}) that is possessed of many descriptions.  One description
in \cite{CFP} uses pairs of finite binary trees.  Parenthesized
expressions are captured by trees.  Thus \(E_1 =
X_1\otimes(X_2\otimes X_3)\) is captured by the tree \(\xy
(0,0);(2,2)**@{-};(6,-2)**@{-}; (4,0);(2,-2)**@{-} \endxy\) and
\(E_2 = (X_1\otimes X_2)\otimes X_3\) is captured by \(\xy
(0,-2);(4,2)**@{-}; (6,0)**@{-}; (2,0);(4,-2)**@{-} \endxy\).  We
can summarize the fact that \(\alpha_0\) has \(E_1\) as source and
\(E_2\) as target by writing \[\alpha_0=\left(\,\,\xy
(0,0);(2,2)**@{-};(6,-2)**@{-}; (4,0);(2,-2)**@{-} \endxy\,\,,\,\,
\xy (0,-2);(4,2)**@{-}; (6,0)**@{-}; (2,0);(4,-2)**@{-}
\endxy\,\,\right).\]

Our isomorphisms connect paired expressions having the same number
of variables, and these paired expressions correspond to pairs of
trees that have the same number of leaves.  The elements of the
group \(F\) are the equivalence classes of all pairs of finite
binary trees in which the two trees in the pair have the same number
of leaves.  It is easiest to explain the equivalence relation put on
such pairs of trees by saying that two pairs are equivalent if they
correspond to instances of the same isomorphism.  The multiplication
of pairs is defined by writing \((T_1, T_2)(T_2,T_3)=(T_1,T_3)\).
This multiplies elements in the reverse order that we have composed
isomorphisms, so the discussion that follows will have some flips in
it.  The arguments that all equivalence classes can be multiplied in
a well defined manner and that the resulting multiplication gives a
group can be found in \cite{CFP}.

Those familiar with Thompson's group \(F\) will recognize the proof
of Theorem \ref{MainThm} as the bulk of yet another proof that \(F\)
has a certain presentation.  With our right-to-left convention for
composing isomorphisms, we end up with the non-standard version of
the presentation that reads \[F=\langle \alpha_0, \alpha_1, \ldots
\mid \alpha_i\alpha_j = \alpha_{j+1}\alpha_i, \quad
\mathrm{whenever} \quad i<j\rangle.\] The usual presentation would
have the relations read \(\alpha_j\alpha_i=\alpha_i\alpha_{j+1}\)
when \(i<j\).

That Thompson's group \(F\) is closely associated to associativity
is well known.  See \cite{dehornoy:assoc} and the end comments of
\cite{mck+thomp}.  Further, given a category \(\scr{C}\) with
multiplication \(\otimes\) and associativity isomorphism \(\alpha\),
it is possible to define a group \(G(\scr{C}, \otimes, \alpha)\)
that will be isomorphic to \(F\) if and only if \((\scr{C}, \otimes,
\alpha)\) is coherent in the sense of \cite{MacLane:coherence} (and
not in our more restrictive sense).  This statement is nothing more
than a checking of definitions.  There is a similar statement
connecting the symmetric, monoidal categories (which combine
associativity and commutativity) with another of Thompson's groups
known as \(V\).  Again, this is just a check of definitions and
repeats a well known connection between \(V\) and the pair
consisting of associativity and commutativity (see
\cite{dehornoy:geom-presentations}).  There is a less trivial
connection between the braided tensor categories of
\cite{joyal+street} and a braided version of \(V\) constructed in
\cite{brin:bv} and \cite{dehornoy:geom-presentations}.  This
connection will be explored elsewhere.

The usual theorem on coherence of associativity, Theorem 3.1 of
\cite{MacLane:coherence}, involves the full power of the operation
\(\otimes\) on natural isomorphisms.  Although not apparent in
\cite{MacLane:coherence}, the main affect of this is to introduce
the left stabilizations of \(\alpha\).  A glance at the pentagon
diagram (3.5) of \cite{MacLane:coherence} shows that the diagram can
be used to express the left stabilization
\(\mathbf{1}\otimes\alpha\) in terms of the right stabilization
\(\alpha\otimes \mathbf{1}\) and instances of \(\alpha\).

As a final remark, we point out that \(F\) has a presentation with
only the generators \(\alpha_0\) and \(\alpha_1\) and only two
relations.  That \(\alpha_0\) and \(\alpha_1\) suffice to generate
follows from the relations given and also from the commutativity of
the diagram \tref{SmallRec}.  If desired, Theorem \ref{MainThm} can
be restated to end with the words: \dots {\itshape unique natural
isomorphism constructed as a composition of instances of
\(\alpha_0\), \(\alpha_1\) and their inverses.}  There is nothing to
be learned from the small number of relations since naturality gives
all the relations that are needed and more.


\begin{thebibliography}{1}

\bibitem{brin:bv}
Matthew~G. Brin, \emph{The algebra of strand splitting. {I}. {T}he algebraic
  structure of the braided {T}hompson group.}, preprint, Binghamton University,
  2004.

\bibitem{CFP}
J.~W. Cannon, W.~J. Floyd, and W.~R. Parry, \emph{Introductory notes on
  {R}ichard {T}hompson's groups}, Enseign. Math. (2) \textbf{42} (1996),
  no.~3-4, 215--256. \MR{98g:20058}

\bibitem{dehornoy:assoc}
Patrick Dehornoy, \emph{The structure group for the associativity identity}, J.
  Pure Appl. Algebra \textbf{111} (1996), no.~1-3, 59--82. \MR{97d:55033}

\bibitem{dehornoy:geom-presentations}
Patrick Dehornoy, \emph{Geometric presentations for {T}hompson's groups and
  related groups}, preprint, University of Caen, 2004.

\bibitem{joyal+street}
Andr{\'e} Joyal and Ross Street, \emph{Braided tensor categories}, Adv. Math.
  \textbf{102} (1993), no.~1, 20--78. \MR{94m:18008}

\bibitem{MacLane:coherence}
Saunders MacLane, \emph{Natural associativity and commutativity}, Rice Univ.
  Studies \textbf{49} (1963), no.~4, 28--46.

\bibitem{mck+thomp}
Ralph McKenzie and Richard~J. Thompson, \emph{An elmentary construction of
  unsolvable word problems in group theory}, Word Problems (Boone, Cannonito,
  and Lyndon, eds.), North-{H}olland, 1973, pp.~457--478.

\end{thebibliography}

\providecommand{\bysame}{\leavevmode\hbox to3em{\hrulefill}\thinspace}
\providecommand{\MR}{\relax\ifhmode\unskip\space\fi MR }
\providecommand{\MRhref}[2]{%
  \href{http://www.ams.org/mathscinet-getitem?mr=#1}{#2}
}
\providecommand{\href}[2]{#2}

\noindent Department of Mathematical Sciences

\noindent State University of New York at Binghamton

\noindent Binghamton, NY 13902-6000

\noindent USA

\noindent email: matt@math.binghamton.edu

\end{document}